%% file: iclr2023_conference.tex
\documentclass{article} 
\usepackage{iclr2023_conference,times}

\input{math_commands.tex}

\usepackage[utf8]{inputenc} 
\usepackage[T1]{fontenc}    
\usepackage{hyperref}       
\usepackage{url}            
\usepackage{booktabs}       
\usepackage{amsfonts}       
\usepackage{nicefrac}       
\usepackage{microtype}

\usepackage{hyperref}
\usepackage{url}
\usepackage{bm}
\usepackage{amsmath}
\usepackage{amsthm}
\usepackage{amssymb}
\usepackage{lipsum}
\usepackage{pifont}
\usepackage{graphicx}
\usepackage{subcaption}
\usepackage{caption}
\usepackage{varwidth}
\usepackage{centernot}
\usepackage{todonotes}
\usepackage{algorithmic, algorithm}
\usepackage{multicol}
\usepackage{multicol}
\usepackage{lipsum}
\usepackage[font=small,labelfont=bf]{caption}
\usepackage[normalem]{ulem} 
\usepackage{float}

\usepackage{hyperref}       
\usepackage{url}            
\usepackage{booktabs}       
\usepackage{amsfonts}       
\usepackage{nicefrac}       
\usepackage{microtype}      
\usepackage{xcolor}         
\usepackage{bm}
\usepackage{amsmath}
\usepackage{amsthm}
\usepackage{amssymb}
\usepackage{lipsum}
\usepackage{pifont}
\usepackage{graphicx}
\usepackage{subcaption}
\usepackage{caption}
\usepackage{varwidth}
\usepackage{centernot}
\usepackage{todonotes}
\usepackage{algorithmic, algorithm}
\usepackage{multicol}
\usepackage{lipsum}

\newtheorem{definition}{Definition}
\newtheorem{remark}{Remark}

\newtheorem{corollary}{Corollary}
\newtheorem{theorem}{Theorem}
\newtheorem{proposition}{Proposition}

\newtheorem{claim}{Claim}
\usepackage{enumerate}
\usepackage{fancyhdr}
\usepackage{datetime}
\usepackage{enumitem}

\title{\Large Gradient Properties of Hard Thresholding Operator}

\iclrfinalcopy

\author{Saeed Damadi, Jinglai Shen \\
  Department of Mathematics and Statistics
  \\
  University of Maryland, Baltimore County\\
  Baltimore, MD 21250 \\
  \texttt{sdamadi1@umbc.edu, shenj@umbc.edu} \\
}

\begin{document}

\maketitle

\begin{abstract}
Sparse optimization receives increasing attention in many applications such as compressed sensing, variable selection in regression problems, and recently neural network compression in machine learning. For example, the problem of compressing a neural network is a bi-level, stochastic, and nonconvex problem that can be cast into a sparse optimization problem. Hence, developing efficient methods for sparse optimization plays a critical role in applications. The goal of this paper is to develop analytical techniques for general, large size sparse optimization problems using the hard thresholding operator. To this end, we study the iterative hard thresholding (IHT) algorithm, which has been extensively studied in the literature because it is scalable, fast, and easily implementable. In spite of extensive research on the IHT scheme, we develop several new techniques that not only recover many known results but also lead to new results. Specifically, we first establish a new and critical gradient descent property of the hard thresholding (HT) operator. Our gradient descent result can be related to the distance between points that are sparse. However, the distance between sparse points cannot provide any information about the gradient in the sparse setting. To the best of our knowledge, the other way around (the gradient to the distance) has not been shown so far in the literature. Also, our gradient descent property allows one to study the IHT when the stepsize is less than or equal to 1/L, where L>0 is the Lipschitz constant of the gradient of an objective function. Note that the existing techniques in the literature can only handle the case when the stepsize is strictly less than 1/L. By exploiting this we introduce and study HT-stable and HT-unstable stationary points and show no matter how close an initialization is to a HT-unstable stationary point (saddle point in sparse sense), the IHT sequence leaves it. Finally, we show that no matter what sparse initial point is selected, the IHT sequence converges if the function values at HT-stable stationary points are distinct, where the last condition is a new assumption that has not been found in the literature. We provide a video of 4000 independent runs where the IHT algorithm is initialized very close to a HT-unstable stationary point and show the sequences escape them. 
\end{abstract}

\section{Introduction}

Solving sparse problems has gained increasing attention in the fields of statistics, finance, and engineering. These problems emerge in statistics as variable selection in linear regression problems \cite{fan2001variable,zou2005regularization,chun2010sparse, desboulets2018review}, mixed-integer programs \cite{bourguignon2015exact,liu2017dual, dedieu2021learning},
portfolio optimization in finance \cite{brodie2009sparse,chang2000heuristics}, compressed sensing in signal processing \cite{foucart2mathematical,eldar2012compressed}, and compressing deep neural networks in machine learning \cite{damadi2022amenable, molchanov2017variational,gale2019state}, just to name a few.
Due to the use of $\ell_0$-(pseudo) norm\footnote{$\ell_0$ is not mathematically a norm because for any norm $\|\cdot \|$ and $\alpha \in \mathbb{R}$, $\|\alpha \bm{\theta} \| = \vert \alpha \vert \|\bm{\theta}\|$, while $\|\alpha \bm{\theta} \|_0 = \vert \alpha \vert \|\bm{\theta}\|_0$ if and only if $\vert \alpha \vert = 1$}, these problems are discontinuous and nonconvex. The $\ell_0$-norm case have been addressed by the hard thresholding (HT) techniques specially the iterative HT (IHT) scheme \cite{blumensath2008iterative,beck2013sparsity,lu2014iterative,zhou2021global}. The Lasso-type, Basic Pursuit(BP)-type, and BP denoising(BPDN)-type problems consider $\ell_1$-norm as a convex approximation of $\ell_0$-norm \cite{tibshirani1996regression,mousavi2019solution}.  Nonconvex approximation of $\ell_0$-norm as $\ell_p$-(pseudo) norm ($0<p<1$) has also been studied well \cite{chartrand2007exact,foucart2009sparsest, lai2011unconstrained,wang2011performance,zheng2017does,won2022unified}. Sparse optimization problems can also be formulated as mixed-integer programs \cite{burdakov2016mathematical}.
Intrinsic combinatorics involved in sparse optimization problems makes it an NP-hard problem (even for a quadratic loss \cite{davis1994adaptive,natarajan1995sparse}) so it is difficult to find a global minimizer. However, greedy algorithms have developed to find local minimizers. To this end, following the ideas of matching pursuit (MP) and orthogonal MP (OMP) \cite{mallat1993matching,pati1993orthogonal} as greedy algorithms, numerous other greedy algorithms have been developed such as stagewise OMP (StOMP) \cite{donoho2012sparse}, regularized OMP (ROMP) \cite{needell2009uniform,needell2010signal}, Compressive Sampling MP (CoSaMP) \cite{needell2009cosamp}, and Gradient Support Pursuit (GraSP) \cite{bahmani2013greedy}. It should be noted that sparse optimization is not restricted to finding a sparse vector.
For example \cite{fornasier2011low,haeffele2014structured,davenport2016overview},
finding a low-rank matrix is considered. The problem of finding a low-rank matrix is a counterpart to finding a sparse vector when it comes to applications dealing with matrices. 
In addition to devising algorithms for solving sparse optimization problems, developing first and second order optimality conditions have also been addressed well \cite{pan2017optimality, bauschke2014restricted, beck2016minimization,li2015first, Lu2015OptimizationOS, pan2015solutions, bucher2018second}. 

The general sparse optimization problem is the following:
\begin{equation}\label{eq:generaloptimization}
\quad
\begin{array}{l}
\min f(\mathbf{x}) \\
\text{s.t. }
C_s \cap \mathcal{X}
\end{array}
\end{equation}
where 
$C_s=\{\mathbf{x} \in \mathbb{R}^n \mid \|\mathbf{x}\|_0 \leq s\}$ (sparsity constraint) is the union of finitely many subspaces of dimension $s$ such that $1 \leq s<n$, $\mathcal{X}$ is a constraint set in $\mathbb{R}^n$, and the objective function $f: \mathbb{R}^n \to \mathbb{R}$ is lower bounded and continuously differentiable, i.e., $C^1$.
In this paper we address a special case of Problem (\ref{eq:generaloptimization}) where $\mathcal{X}=\mathbb{R}^n$ as follows:
\begin{equation}\label{eq:optimizationproblem}
(\text{P}):  
\quad
\begin{array}{l}
\min f(\mathbf{x}) \\
\text{s.t. }
\mathbf{x} \in C_s
\end{array}
\end{equation}
To address Problem (\ref{eq:optimizationproblem})  the following fundamental questions arise:
\begin{enumerate}[label=(Q\arabic*)]
\item What are the necessary/sufficient conditions for a local/global minimizer of Problem (\ref{eq:optimizationproblem})?
\item What are the characteristics of accumulation points of algorithms solving Problem (\ref{eq:optimizationproblem})?
\item 
Under what condition(s) does an accumulation point become a local/global minimizer?
\item 
If an accumulation point is a local/global minimizer, what is the rate of convergence?
\end{enumerate}
\input{iht}
By considering the IHT algorithm, we will answer the above questions.
This algorithm has been extensively studied in the literature. It was originally devised for solving compressed sensing problems in 2008 \cite{blumensath2008iterative,blumensath2009iterative}. Since then, there has been a large body of literature studying the IHT-type algorithms from different standpoints. For example, \cite{beck2013sparsity,lu2014iterative, Lu2015OptimizationOS,pan2017convergent,zhou2021global} consider convergence of iterations, \cite{jain2014iterative, liu2020between} study the limit of the objective function value sequence, \cite{liu2017dual,zhu2018lagrange} address duality, \cite{zhou2020subspace, zhao2021lagrange} extend it to Newton's-type IHT, 
\cite{chen2016accelerated,li2016nonconvex,liang2020effective,zhou2018efficient} consider the stochastic version, \cite{blumensath2012accelerated,khanna2018iht,vu2019accelerating,wu2020accelerated} address accelerated IHT, and \cite{wang2019fast, bahmani2013greedy} solve logistic regression problem using the IHT. 
\subsection*{Summary of Contributions}
%
By considering the IHT Algorithm \ref{alg:IHT} for Problem (\ref{eq:optimizationproblem}), we develop the following results:
\begin{itemize}
\item
We establish a new critical gradient descent property of the hard
thresholding (HT) operator that has not been found in the literature. Our gradient descent result can be related to the distance between points that are sparse. However, the distance between sparse points cannot provide any information about the gradient in the sparse setting. To the best of our knowledge, the other way around (the gradient to the distance) has not been shown so far in the literature. This property allows one to study the IHT when the stepsize is less
than or equal to $1/L$, where $L > 0$ is the Lipschitz constant of the gradient of an objective function. Note that the existing techniques in the literature can only handle the case when the stepsize is
strictly less than $1/L$. As an example, one can refer to \cite{liu2020between} that needs the stepsize to be greater than or equal to $1/L$.  
\item
We introduce the notion of \textit{HT-stable/unstable stationary} points.
Using them we establish the \textit{escapability} property of \textit{HT-unstable stationary} points (saddle point the sparse sense) and local \textit{reachability} property of strictly \textit{HT-stable stationary} points. We provide a video of 4000 independent runs where the IHT algorithm is initialized very close to a HT-unstable stationary point and show the sequences escape them. 
\item
We also show that the IHT sequence converges globally under a new assumption that has not been found in the literature. In addition, Q-linearly convergence of the IHT algorithm towards a local minimum when the objective function is both RSS and restricted strictly convex is shown.
\end{itemize}
According to our results, we address (Q1) and (Q2) by  establishing a new gradient descent property of the hard thresholding (HT) operator and introducing 
the notion of \textit{HT-stable/unstable stationary} points. By considering RSS,  restricted strictly convex, and RSC properties we address (Q3) and (Q4). Table \ref{tab:comparison} is provided to compare our results with those in the literature. It shows what has been done chronologically and demonstrates our results.
\section{Related work}

To answer (Q1) \cite{beck2013sparsity} introduces \textit{$L$-stationarity} property as a necessary condition for an optimal solution of Problem (\ref{eq:optimizationproblem}). The \textit{$L$-stationarity} property is defined when the gradient of the objective function is Lipschitz. Also, \cite{beck2013sparsity} addresses (Q2) by showing any accumulation point of the IHT algorithm is $L$-stationary.  
Lu in \cite{lu2014iterative} restricts the objective function to be convex and shows that the IHT sequence converges to a local minimum when the objective function is regularized by $\ell_0$-norm and $\mathcal{X}$ is a box constraint. Jain et al. \cite{jain2014iterative} put more restriction on the objective value function and show that
the objective value function sequence generated by the IHT algorithm converges to a value attained under a more restricted sparsity constraint. The restrictions used in \cite{jain2014iterative} are Restricted Strong Smoothness (RSS) and Restricted Strong Convexity (RSC). The RSS and RSC properties are introduced by \cite{negahban2012unified} and first used by \cite{bahmani2013greedy} for sparsity optimization problems. Currently, they have become standard restrictions for analyzing sparsity optimization problems. Under RSS and RSC properties for the objective function, one is able to address (Q3) and (Q4).

Finding a closed-form expression for $P_{C_s \cap \mathcal{X}}$ when $\mathcal{X}$ is an arbitrary set is difficult. However, \cite{beck2016minimization}
shows orthogonal projection of a point onto $C_s \cap \mathcal{X}$ can be efficiently computed
when $\mathcal{X}$ is a symmetric closed convex set. In this context, two types of sets are of interest: nonnegative symmetric sets and sign free sets. To address (Q1) in a more generalized setting, 
Beck and Hallak \cite{beck2016minimization} characterize \textit{$L$-stationary} points of Problem (\ref{eq:generaloptimization}) when $\mathcal{X}$ is either nonnegative symmetric set or sign free. Also, Lu in \cite{Lu2015OptimizationOS} considers the same setting as \cite{beck2016minimization} and introduces a new optimality condition that is stronger than \textit{$L$-stationary}. He devises a Nonmonotone Projected Gradient (NPG) algorithm and shows an accumulation of the NPG sequence is the global optimal of Problem (\ref{eq:generaloptimization}). Pan et al. \cite{pan2017optimality} consider Problem (\ref{eq:generaloptimization}) when $\mathcal{X}=\mathbb{R}^n_{+}$. They develop an Improved IHT algorithm (IIHT) that employs the Armijo-type stepsize rule. They show when the objective function is RSS and RSC, the IIHT sequence converges to a local minimum. A recent work by Zhou et al. \cite{zhou2021global} develops Newton Hard-Thresholding Pursuit (NHTP) for solving problem (\ref{eq:optimizationproblem}). They show that when accumulation points of the NHTP sequence are \textit{$L$-stationary} and are isolated, the sequence converges with a locally Q-quadratic rate. Table \ref{tab:comparison} compares current results in the literature.
\input{comparison}
\section{Definitions}
We provide some definitions that will be used throughout the paper. These definitions are the HT operator (HTO) and HTO inequality, RSS and RSC functions.
%
\begin{definition}[Restricted Strong Smoothness (RSS)]\label{def:rss}
A differentiable function $f: \mathbb{R}^n \to \mathbb{R}$ is said to be restricted strongly smooth with modulus $L_s>0$ or is $L_s$-RSS if
\begin{equation}\label{eq:rss}
f(\mathbf{y}) \leq f(\mathbf{x}) + \langle \nabla f(\mathbf{x}) , \mathbf{y}-\mathbf{x} \rangle + \frac{L_{s}}{2}\|\mathbf{y}-\mathbf{x}\|_2^2 \quad \forall \mathbf{x},\mathbf{y} \in \mathbb{R}^n \text{ such that } \|\mathbf{x}\|_0 \leq s,\|\mathbf{y}\|_0\leq s.
\end{equation}
\end{definition}
\begin{definition}[Restricted Strong Convexity (RSC)]\label{def:rsc}
A differentiable function $f: \mathbb{R}^n \to \mathbb{R}$ is said to be restricted strongly convex with modulus $\beta_s>0$ or is $\beta_s$-RSC if
\begin{equation}\label{eq:rsc}
f(\mathbf{y}) 
\geq
f(\mathbf{x}) + \langle \nabla f(\mathbf{x}) , \mathbf{y}-\mathbf{x} \rangle + \frac{\beta_s}{2}\|\mathbf{y}-\mathbf{x}\|_2^2 
\quad 
\forall \mathbf{x},\mathbf{y} \in \mathbb{R}^n
\text{ such that } ||\mathbf{x}||_0 \leq s,||\mathbf{y}||_0\leq s.
\end{equation}
\end{definition}

\begin{definition}
[The HT operator]
\label{def:hardthresholding}
The HT operator $H_s(\cdot)$ denotes the orthogonal projection onto multiple subspaces of $\mathbb{R}^n$ with dimension $1 \leq s<n$, that is,
\begin{equation}\label{eq:hardthreshold}
    H_s(\mathbf{x}) \in \arg\min_{\|\mathbf{z}\|_0\leq s }\|\mathbf{z}-\mathbf{x}\|_2.
\end{equation}
\end{definition}
\begin{claim}\label{claim:tops}
The HT operator keeps the $s$ largest entries of its input in absolute values.
\end{claim}
For a vector $\mathbf{x} \in \mathbb{R}^n$, $\mathcal{I}^{\mathbf{x}}_s \subset \{1,\dots, n\}$ denotes the set of indices corresponding to the first $s$ largest elements of $\mathbf{x}$ in absolute values. For example $H_2([1,-3,1]^{\top})$ is either $[0,-3,1]^{\top}$ or $[1,-3,0]^{\top}$ where $\mathcal{I}^{\mathbf{y}}_2=\{2,3\}$ and $\mathcal{I}^{\mathbf{y}}_2=\{1,2\}$, respectively. Therefore, the output of it may not be unique. This clearly shows why HTO is not a convex operator and why there is an inclusion in (\ref{eq:hardthreshold}) not an inequality.
\section{Results}
We consider solving Problem (\ref{eq:optimizationproblem}) using the IHT Algorithm \ref{alg:IHT} and develop results on the HT operator. Using them, the behavior of the IHT sequence generated by Algorithm \ref{alg:IHT} is characterized. 
Towards this end, statements of the main results are provided and all the technical proofs are postponed to the Appendix for the reviewers.

\subsection{Gradient descent property}

First, we establish a new and critical gradient descent property of the hard thresholding
(HT) operator.
\begin{theorem}\label{theorem:deterministicunion}
Let $f: \mathbb{R}^n \rightarrow \mathbb{R}$ be a differentiable function that is $L_{s}$-RSS, $\mathbf{y} \in H_s(\mathbf{x}-\gamma \nabla f(\mathbf{x}))$ with any $\mathcal{I}_s^{\mathbf{y}}$ and $0<\gamma\le\frac{1}{L_s}$, and $\mathbf{x}$ be a sparse vector such that $\|\mathbf{x}\|_0\leq s$ with any $\mathcal{I}_s^{\mathbf{x}}$. Then,
\begin{equation}\label{eq:deterministicunion}
\frac{\gamma}{2}(1-L_{s}\gamma)
\|\nabla_{\mathcal{I}_s^{\mathbf{x}} \cup \mathcal{I}_s^{\mathbf{y}}}f(\mathbf{x})\|_2^2  \leq f(\mathbf{x}) - f(\mathbf{y})  
\end{equation}
where $\nabla_{\mathcal{I}_s^{\mathbf{x}} \cup \mathcal{I}_s^{\mathbf{y}}}f(\mathbf{x})$ is the restriction of the gradient vector to the union of the index sets $\mathcal{I}_s^{\mathbf{x}}$ and $ \mathcal{I}_s^{\mathbf{y}}$. 
\end{theorem}

Theorem \ref{theorem:deterministicunion} provides a lower bound on the difference between the current function value evaluated at $\mathbf{x}$ and the one evaluated at the updated point provided by the HTO, i.e., $\mathbf{y}$. Note that, $\mathbf{y}$ may not be a unique vector that has $s$ nonzero elements. Nonetheless, as stated in Theorem \ref{theorem:deterministicunion}, Inequality (\ref{theorem:deterministicunion}) holds for any $\mathbf{y}$ that might be the output of the HTO. As one clearly see, the descent can only be characterized by looking at the entries of the gradient that are restricted to the union of the $s$ largest elements in both $\mathbf{x}$ and $\mathbf{y}$. The rest of the gradient can be ignored. Since one may be interested in characterizing the descent using the distance between $\mathbf{x}$ and $\mathbf{y}$, we provide the following corollary.

\begin{corollary}\label{cor:neiboringerror}
Assume all the assumptions in Theorem \ref{theorem:deterministicunion} hold, then,
\begin{equation}\label{eq:neighboringerror}
\frac{1-L_{s}\gamma}{6\gamma}
\|\mathbf{y}-\mathbf{x}\|_2^2
\leq
\frac{\gamma}{2}(1-L_{s}\gamma)
\|\nabla_{\mathcal{I}_s^{\mathbf{x}} \cup \mathcal{I}_s^{\mathbf{y}}}f(\mathbf{x})\|_2^2  \leq f(\mathbf{x}) - f(\mathbf{y}) 
\end{equation}
\end{corollary}

The above result shows the superiority of our gradient result because our gradient result can be related to the distance of points that are sparse. However, the distance between sparse points cannot provide any information about the gradient. To the best of our knowledge, the other way around (the gradient to the distance) has not been shown so far in the literature.

Algebraically speaking, characterizing a descent of the function value solely with the information of the current iterate, i.e., $\mathbf{x}$, is of more interest. To this end, we provide another corollary to Theorem \ref{theorem:deterministicunion} that ties the descent to $\mathbf{x}$ only.
\begin{corollary}\label{cor:deterministicx}
Assume all the assumptions in Theorem \ref{theorem:deterministicunion} hold. Then, the norm of the gradient restricted to any $\mathcal{I}_s^{\mathbf{x}}$ can be bounded as follows:
\begin{equation}\label{eq:deterministicgammaeqL}
\frac{\gamma}{2}
\|\nabla_{\mathcal{I}_s^{\mathbf{x}}}f(\mathbf{x})\|_2^2  \leq f(\mathbf{x}) - f(\mathbf{y})  
\end{equation}
\end{corollary}

By this point, we have shown that applying the HTO once, can result in smaller function value provided the gradient over the $s$ largest entries of $\mathbf{x}$ are nonzero. This can be utilized to show the sequence generated by the IHT algorithm is nonincreasing. Specially, if the generated sequence has an accumulation point, the objective value function sequence converges to the objective value of the accumulation point.
\footnote{A sequence may not converge but it may have an accumulation point. For example $1,-1,1,-1,\dots$ is not a convergent sequence but it has two accumulation points.}

\begin{corollary}\label{cor:fconvergence}
Let $f: \mathbb{R}^n \rightarrow \mathbb{R}$ be a bounded below differential function that is $L_{s}$-RSS and $\big(\mathbf{x}^k\big)_{k\geq 0}$ be the IHT sequence $\big(\mathbf{x}^k\big)_{k\geq 0}$
with
$
0
<
\gamma
\leq
\frac{1}{L_s}$.
Then, $\Big(f(\mathbf{x}^{k})\Big)_{k\geq 0}$ is nonincreasing and converges. Also, if
$\mathbf{x}^*$ is an accumulation point of $\big(\mathbf{x}^{k}\big)_{k\geq 0}$
then $\Big(f(\mathbf{x}^{k})\Big)_{k\geq 0}\to f(\mathbf{x}^*)$.
\end{corollary}

Next, we look at \textit{basic stationary points} of Problem (\ref{eq:optimizationproblem}) and show their properties.

\subsection{Optimality condition based on the HT properties}\label{subsub:IHTstability}
In this subsection, we will show that not all \textit{basic stationary} points of Problem (\ref{eq:optimizationproblem}) are reachable when the IHT algorithm is run. To do so, the notion of \textit{HT stationary} points are introduced as follows.
\subsubsection{\textit{HT stationary} points}
\begin{definition}\label{def:HTstable}
For a given constant $\gamma>0$, we say that a sparse vector $\mathbf{x}^* \in C_s$ is \textit{HT-stable} stationary point of Problem (\ref{eq:optimizationproblem}) associated with $\gamma$ if $\nabla_{\text{supp}(\mathbf{x}^*)} f(\mathbf{x}^*)=0$, and 

\begin{equation}\label{eq:HTstable}
\min
\Big(
|x^*_i|:
i\in \mathcal{I}^{{\mathbf{x}^*}}_s
\Big)
\ge
\gamma
\max\Big(
|
\nabla_j f(\mathbf{x}^*)
|:
j \notin \text{supp}(\mathbf{x}^*
)\Big)
=
\gamma
\|
\nabla_{({\text{supp}(\mathbf{x}^*
)})^c} f(\mathbf{x}^*)
\|_{\infty}.
\end{equation}

(Note that $\min
\Big(
|x^*_i|:
i\in \mathcal{I}^{\mathbf{x}^*}_s
\Big)$ is unique and does not depend on the choice $\mathcal{I}^{\mathbf{x}^*}_s$.) If $\nabla_{\text{supp}(\mathbf{x}^*)} f(\mathbf{x}^*)=0$ but (\ref{eq:HTstable}) fails, we say that $\mathbf{x}^*$ is \textit{HT-unstable} stationary point with $\gamma$. Moreover, if $\nabla_{\text{supp}(\mathbf{x}^*)} f(\mathbf{x}^*)=0$ and (\ref{eq:HTstable}) holds strictly, namely,

\begin{equation}\label{eq:HTstrictly}
\min
\Big(
|x^*_i|:
i\in \mathcal{I}^{\mathbf{x}^*}_s
\Big)
>
\gamma
\max\Big(
|
\nabla_j f(\mathbf{x}^*)
|:
j \notin \text{supp}(\mathbf{x}^*
)\Big)
\end{equation}
we say that $\mathbf{x}^*$ is a strictly \textit{HT-stable} stationary point associated with $\gamma$.
\end{definition}
Note that when $\|\mathbf{x}^*\|_0=s$, $\mathcal{I}^{\mathbf{x}^*}_s$ is unique and equals $\text{supp}(\mathbf{x}^*)$ such that $\text{supp}(\mathbf{x}^*)$ in the above definition can be replaced by  $\mathcal{I}^{\mathbf{x}^*}_s$. Moreover, if $\mathbf{x}^*$ is a strictly \textit{HT-stable} stationary point, then we must have $\mathcal{I}^{\mathbf{x}^*}_s=\text{supp}(\mathbf{x}^*)$ (or equivalently $\|\mathbf{x}^*\|_0=s$) because otherwise, $0=\min
\Big(
|x^*_i|:
i\in \mathcal{I}^{{\mathbf{x}^*}}_s
\Big)
>
\gamma
\|
\nabla_{({\text{supp}(\mathbf{x}^*
)})^c} f(\mathbf{x}^*)
\|_{\infty}$
which is impossible.
\begin{remark}
As stated in the Definition \ref{def:HTstable}, a basic stationary point is a point whose gradient is zero over the nonzero elements. For example, suppose $\tilde{\mathbf{x}}=[0, 4, 0, 2]^{\top} \in \mathbb{R}^4$ is a basic stationary point. Then $\nabla f(\tilde{\mathbf{x}})=[c_1, 0, c_2, 0]^{\top}$ where $c_1,c_2$ are scalars. The main idea of the HT-stable stationary point is that it has to be a basic stationary point. In other words $\tilde{\mathbf{x}}$ can be a basic stationary point but not a HT-stable stationary point. This is the analogue of the non-sparse optimization where a point $\hat{\mathbf{x}}$ whose gradient is zero, i.e., $\nabla f(\hat{\mathbf{x}})=0$ may not be necessary a local or global minimizer. It can be a saddle point.
\end{remark}
\begin{remark}
The main message of Definition \ref{def:HTstable} is the following: ``only by looking at the gradient restricted to the nonzero entries of a \textit{basic feasible point}, one cannot say whether it is a local minimizer of Problem (\ref{eq:optimizationproblem}) or not''.
\end{remark}
An \textit{HT-stable} stationary point associated with $\gamma$ is equivalent to the $\frac{1}{\gamma}$-stationary point of Problem (\ref{eq:optimizationproblem}) defined in \cite[Definition 2.3]{beck2013sparsity}. Thus, by \cite[Lemma 2.2]{beck2013sparsity}, $\mathbf{x}^*$ is a \textit{HT-stable} point if and only if $\mathbf{x}^* \in H_s(\mathbf{x}^* -\gamma \nabla f(\mathbf{x}^*))$. The notion of a HT-unstable stationary point is novel and is a key point for proving Theorem 2. Theorem 2 is the foundation for the proof of part b) of Theorem 3 as well as Theorem 4 which characterizes the accumulation point of the IHT sequence. In addition, we have introduced another stationary point, namely strictly HT-stable stationary which is a crucial concept for local convergence of the IHT sequence. 

In the following, we present a result that characterizes a \textit{HT-unstable} stationary point. In essence, the following result shows that there always exists a neighborhood around a sparse \textit{HT-unstable} stationary point whose gradient is zero over the nonzero elements and one can decrease the function value by going towards the direction of any nonzero coordinates.
\begin{theorem}\label{theorem:HTunstable}
Suppose $f:\mathbb{R}^n\to\mathbb{R}$ is $C^1$ and $L_s$-RSS. Given any $0<\gamma \leq \frac{1}{L_s}$, if a vector $\tilde{\mathbf{x}}\in C_s$ is such that $\nabla_{\text{supp}(\tilde{\mathbf{x}})}f(\tilde{\mathbf{x}})=0$ and $\min \big( |\tilde{x}_i|: i \in \mathcal{I}^{\tilde{\mathbf{x}}}_s\big)<\gamma \|\nabla_{(\text{supp}(\tilde{\mathbf{x}}))^c}f(\tilde{\mathbf{x}})\|_{\infty}$ for some $\mathcal{I}^{\tilde{\mathbf{x}}}_s$, then there exist a constant $\nu > 0$ and a neighborhood $\mathcal{N}$ of $\tilde{\mathbf{x}}$ such that $f(\mathbf{y}) \leq f(\mathbf{x})-\nu$ for any $\mathbf{x} \in \mathcal{N}\cap C_s$ and any $\mathbf{y} \in H_s(\mathbf{x}-\gamma \nabla f(\mathbf{x})).$
\end{theorem}

The above result leads to the following necessary optimality conditions for a global minimizer of Problem (\ref{eq:optimizationproblem}) in terms of hard thresholding operator $H_s$. 
For the case where $\gamma=1/L_s$, i.e., part b), one needs to use Theorem 2. Indeed, to the best of our knowledge, no proof has not been found for it in the literature. Essentially, establishing the result in part b) is one of our contributions. For the case $\gamma < 1/L_s$ it is proven that $\mathbf{x}^*=H_s(\mathbf{x}^*
-\gamma \nabla f(\mathbf{x}^*))$. 
Note that the condition for $0 < \gamma < \frac{1}{L_s}$ has been obtained in \cite[Theorem 2.2]{beck2013sparsity} without using gradient properties of the HT operator.
\begin{theorem}\label{theorm:global}
Suppose $f:\mathbb{R}^n\to\mathbb{R}$ is $L_s$-RSS and $\mathbf{x}^*$ is a global minimizer. Then, $\mathbf{x}^*$ is a HT-stable (or $\frac{1}{\gamma}$-) stationary point for any $0< \gamma \leq \frac{1}{L_s}$. Particularly, the following hold:
\begin{enumerate}[label=\alph*]
\item) 
For any $0<\gamma < \frac{1}{L_s}$, 
$\mathbf{x}^*=H_s(\mathbf{x}^*
-\gamma \nabla f(\mathbf{x}^*))$.
\item) 
For $\gamma =\frac{1}{L_s}$, $\mathbf{x}^* \in H_s(\mathbf{x}^*
-\gamma \nabla f(\mathbf{x}^*))$.
\end{enumerate}
\end{theorem}

The following result shows that any accumulation point of an IHT sequence must be a \textit{HT-stable} stationary point.

\begin{theorem}\label{theorm:IHTaccumulation}
Let $f:\mathbb{R}^n\to\mathbb{R}$ be $L_s$-RSS and $C^1$ function. Suppose $f$ is bounded below on $C_s$. Consider an IHT sequence $\big(\mathbf{x}^k\big)_{k\geq 0}$ associated with an arbitrary $\gamma \in (0, \frac{1}{L_s}]$, and let $\mathbf{x}^*$ be an accumulation point of $\big(\mathbf{x}^k\big)_{k\geq 0}$. Then, $\mathbf{x}^*$ is a HT-stable stationary point of Problem (\ref{eq:optimizationproblem}). 
\end{theorem}

\begin{remark}
\normalfont
The above theorem shows that any accumulation point of an IHT sequence is a \textit{HT-stable} stationary point of Problem (\ref{eq:optimizationproblem}). Since each \textit{HT-stable} stationary point must be a basic stationary point, one can observe that any accumulating point $\mathbf{x}^*$ of an IHT sequence must satisfy $\nabla_{\text{supp}(\mathbf{x}^*)} f(\mathbf{x}^*)=0$ when $\|\mathbf{x}^*\|_0=s$, or $\nabla f(\mathbf{x}^*)=0$ when $\|\mathbf{x}^*\|_0<s$.
\end{remark}
The following result pertains to the objective function values of \textit{HT-stable} and \textit{HT-unstable stationary} points.

\begin{corollary}\label{cor:stablevsunstable}
Let $f:\mathbb{R}^n\to\mathbb{R}$ be $L_s$-RSS and $C^1$ function. Suppose that every (nonempty) sub-level set of $f$ contained in $C_s$ is bounded, i.e., for any $\alpha \in \mathbb{R}$, $\{x \in C_s | f(\mathbf{x})\leq \alpha\}$ is bounded (and closed). Consider an arbitrary $\gamma \in (0, \frac{1}{L_s}]$. For any HT-unstable stationary point $\mathbf{x}^*$ associated with $\gamma$, there exists a HT-stable stationary point $\tilde{\mathbf{x}}^*$ associated with $\gamma$ such that $f(\mathbf{x}^*)>f(\tilde{\mathbf{x}}^*)$.
\end{corollary}

Based on the above corollary, it is easy to see that if there are finitely many \textit{HT-unstable} stationary points (happens when the function is RSC), then there is a \textit{HT-stable} stationary point $\tilde{\mathbf{x}}^*$ such that $f(\mathbf{x}^*)>f(\tilde{\mathbf{x}}^*)$ for any \textit{HT-unstable} stationary point $\mathbf{x}^*$.

The following result provides sufficient conditions for the convergence of an IHT sequence. Corollary \ref{cor:IHTconvergence} aims to remove any restrictions on the initial condition. This corollary shows that no matter what initial condition in $C_s$ is selected, the IHT sequence will converge to a HT-stable stationary point. Note that we say a \textit{HT-stable/unstable} stationary point $\mathbf{x}^*$ associated with $\gamma \in (0, \frac{1}{L_s}]$ is isolated if there exists a neighborhood $\mathcal{N}$ of $\mathbf{x}^*$ such that $\mathcal{N}$ does not contain any HT stationary point other than $\mathbf{x}^*$.

\begin{corollary}\label{cor:IHTconvergence}
Let $f:\mathbb{R}^n\to\mathbb{R}$ be $L_s$-RSS and $C^1$ function. Suppose that every (nonempty) sub-level set of $f$ contained in $C_s$ is bounded. Consider an arbitrary $\gamma \in (0,\frac{1}{L_s}]$. Assume that
\begin{enumerate}[label=A.\arabic*]
\item:
For any two distinct HT-stable stationary points $\mathbf{x}^*$ and $\mathbf{y}^*$ associated with $\gamma$, $f(\mathbf{x}^*) \neq f(\mathbf{y}^*)$. Then, for any $\mathbf{x}^0 \in C_s$, the IHT sequence $\big(\mathbf{x}^k\big)_{k\geq 0}$ converges to a HT-stable stationary point associated with $\gamma$. This convergence results also hold under the following assumption:
\begin{enumerate}[label=A.\arabic*]
\setcounter{enumii}{1}
\item:
when $0<\gamma < \frac{1}{L_s}$, each HT-stable stationary point associated with $\gamma$ is isolated.
\end{enumerate}
\end{enumerate}
\end{corollary}

The following corollary shows that any IHT sequence always ``escape" from a \textit{HT-unstable} stationary point.

\begin{corollary}\label{cor:escapability}
Let $f:\mathbb{R}^n\to\mathbb{R}$ be $L_s$-RSS and $C^1$ function. Suppose $f$ is bounded below on $C_s$. For any given $\gamma \in (0, \frac{1}{L_s}]$ and any HT-unstable stationary point $\mathbf{x}^*$ associated with $\gamma$, there exists a neighborhood $\mathcal{N}$ of $\mathbf{x}^*$ such that for any IHT sequence starting from any $\mathbf{x}^0 \in C_s$, there exists $N \in \mathbb{N}$ such that $\mathbf{x}^k \notin \mathcal{N} \cap C_s$ for all $k \geq N$.
\end{corollary}

The next result shows the attraction towards a strictly HT-stable stationary point in a neighborhood of such a stationary point. In what follows, for each index subset $\mathcal{J}$ with $|\mathcal{J}|=s$, a subspace $\mathcal{S_J}:=\{\mathbf{x} \in \mathbb{R}^n \mid \mathbf{x}_{\mathcal{J}^c}\}$ associated with $\mathcal{J}$ is defined.
Clearly, $C_s$ is the union of $\mathcal{S_J}$'s for all $\mathcal{J}$'s with $|\mathcal{J}|=s$.
\begin{proposition}\label{prop:reachability}
Let $f:\mathbb{R}^n\to\mathbb{R}$ be $L_s$-RSS and $C^1$ function. Suppose $f$ is bounded below on $C_s$ and $f$ is strictly convex on $\mathcal{S_J}$ for any index subset $\mathcal{J}$ with $|\mathcal{J}|=s$. Let $\mathbf{x}^*$ be a strictly HT-stable stationary point associated with any given $\gamma \in (0, \frac{1}{L_s}]$. Then there exists a neighborhood $\mathcal{B}$ of $\mathbf{x}^*$ such that for every $\mathbf{x}^0 \in \mathcal{B} \cap C_s$, the IHT sequence $\big(\mathbf{x}^k\big)_{k\geq 0}$ converges to $\mathbf{x}^*$. Moreover, if $f$ is strongly convex on $\mathcal{J}$ for every index subset  $\mathcal{J}$ with $|\mathcal{J}|=s$, then for $\mathbf{x}^0 \in \mathcal{B} \cap C_s$, the IHT sequence $\big(\mathbf{x}^k\big)_{k\geq 0}$ Q-linearly converges to $\mathbf{x}^*$.
\end{proposition}
Next, we provide an example to show the \textit{escapability} property of \textit{HT-unstable} points.
\section{Simulation}
\begin{figure}[!htb]
   \begin{minipage}{0.48\textwidth}
     \centering
     \includegraphics[width=.95\linewidth]{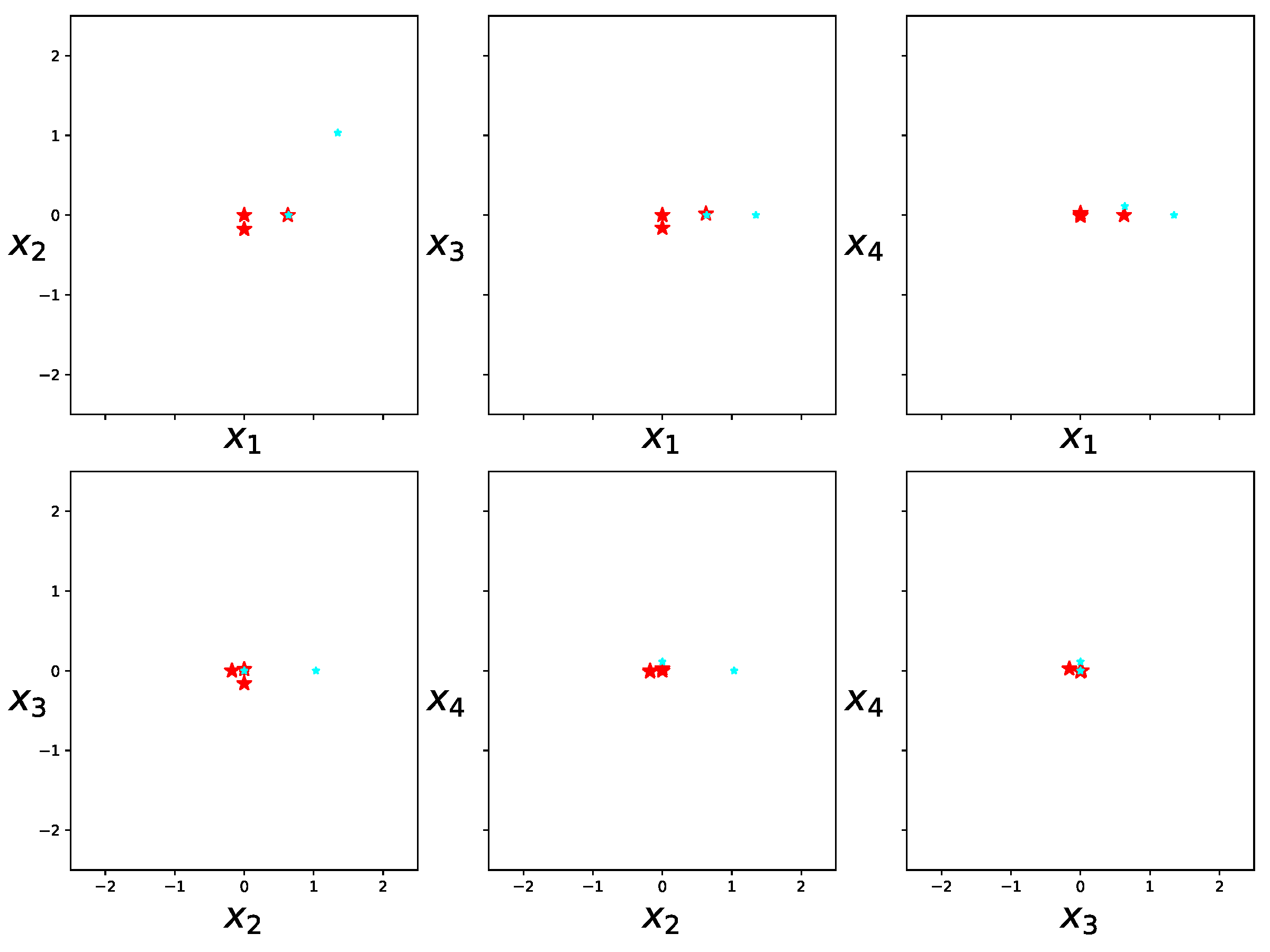}
     \caption{Illustration of  \textit{HT-stationary} points on 2-d plains.}\label{fig:HT-stationary}
   \end{minipage}\hfill
   \begin{minipage}{0.48\textwidth}
     \centering
     \includegraphics[width=.95\linewidth]{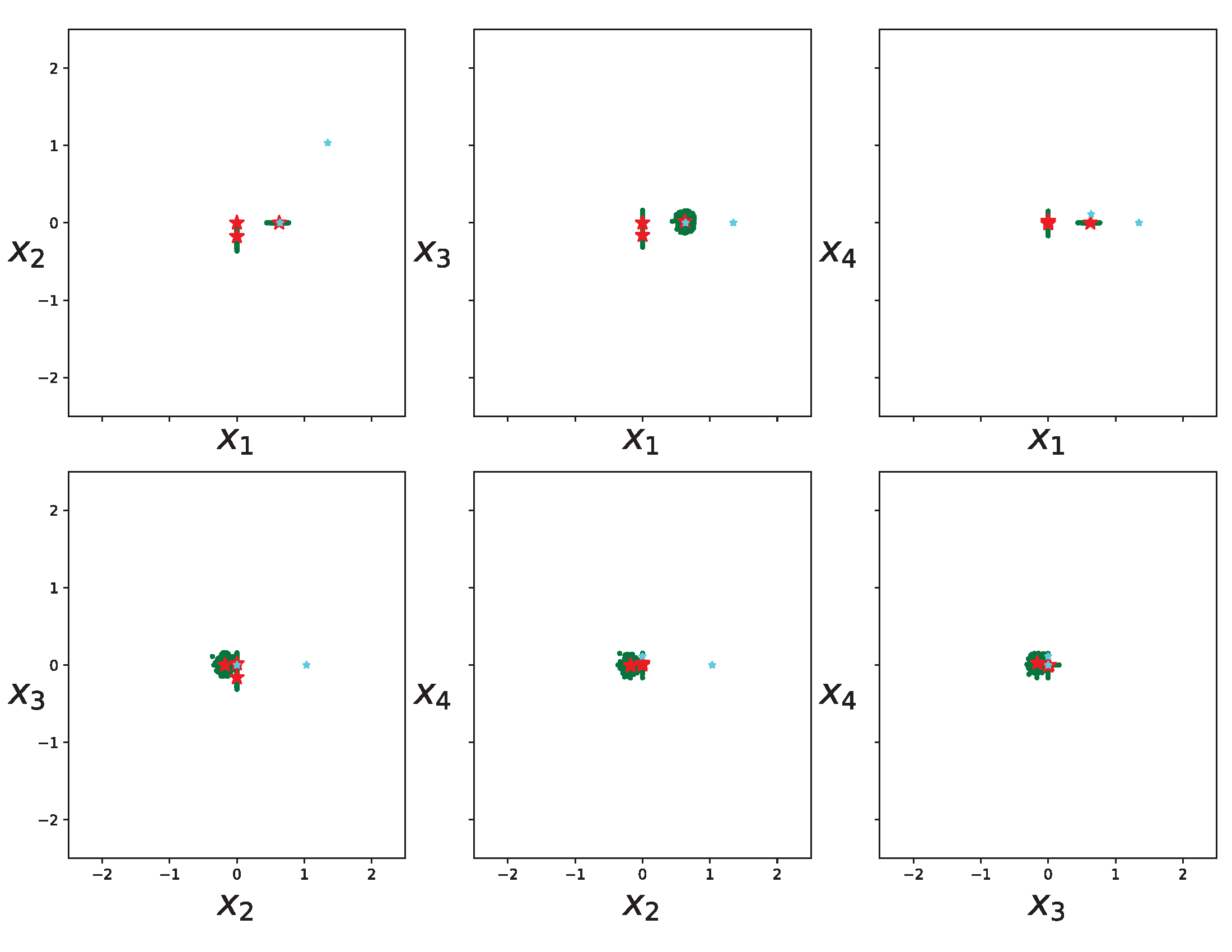}
     \caption{Illustration of  4,000 initialization close to four \textit{HT-unstable stationary} points.}\label{fig:unstable_init}
   \end{minipage}
\end{figure}

To elaborate on theoretical results including 
Corollary \ref{cor:fconvergence}, 
Theorem \ref{theorem:HTunstable}, the notion of \textit{HT-stationary} points, Corollary \ref{cor:escapability}  which shows \textit{escapability} property of \textit{HT-unstable stationary} points, and Proposition \ref{prop:reachability} which shows \textit{Reachability} to \textit{HT-stable stationary} points, we use a quadratic function 
$f(\mathbf{x})=\frac{1}{m}\sum_{i=1}^m(A_{i\bullet}\mathbf{x}-y_i)^2=\frac{1}{m}\|A\mathbf{x}-\mathbf{b}\|^2$
where $\mathbf{A} \in \mathbb{R}^{m\times n}$, $A_{i\bullet}$ is the $i$-th row of $A$, $\mathbf{x} \in \mathbb{R}^n$ is the optimization variable, and $\mathbf{y} \in \mathbb{R}^m$ is the target. This function is both RSS and RSC so both Corollary \ref{cor:escapability} and Proposition \ref{prop:reachability} follow. To better visualize the process, we let $m=n=4$ and $s=2$. Therefore, there are six \textit{HT-stationary} points where the gradient over the nonzero elements is zero. We  use Pytorch \citep{paszke2019pytorch} to select the matrix $A$ and $\mathbf{y}$. By setting the random seed to be 45966 we draw a $4\times 4$ matrix $A$ whose elements are standard normal. Keeping the same seed, we generate $\mathbf{y}$. The following would be $A$ and $\mathbf{y}$:
$$
A=
\begin{bmatrix}
-1.0655 &  0.2249 & -0.0897 &  0.1876 \\
 1.1627 & -1.1229 & -0.0823 & -0.3059 \\
-0.2011 &  0.5342 & -0.0551 & -1.3459 \\
 0.2308 & -0.6404 & -0.7468 &  0.0378
\end{bmatrix}
,\quad
\mathbf{y}=
\begin{bmatrix}
-1.7861\\ -0.3556\\ -0.1881\\  0.3896
\end{bmatrix}
$$
The restricted Lipschitz constant, i.e., $L_s$, for the above quadratic function is $\frac{2}{m}\times \lambda_{max}(A^{\top}A)$ where $\lambda_{max}$ is the maximum eigenvalue of $A^{\top}A$. Thus, for the above choice of $A, \mathbf{y}$, the maximum allowable stepsize is $\gamma=\frac{1}{L_s}=0.06$. Once, $\gamma$ is fixed, one can determine stability of each stationary point.
The following are \textit{HT-stationary} points along with their stability status as well as the gradient status of each \textit{HT-stationary} point.
As you can see, the gradient corresponding to nonzero elements in \textit{HT-stationary} point are zero:
$$
\begin{bmatrix}
No. & x_1 & x_2 & x_3 & x_4 & g_1 & g_2 & g_3 & g_4 & HT-stability\\
1 & 1.3474 & 1.0331 & 0 & 0 & 0 & 0 & 0.2060 & -0.3916 & \text{strictly HT-stable}\\
2 & 0.6278 & 0 & 0.0177 & 0 & 0 & -0.3843 & 0 & -0.1070 & \text{HT-unstable}\\
3 & 0.6387 & 0 & 0 & 0.1123 & 0 & -0.4189 & -0.0029 & 0 & \text{strictly HT-stable}\\
4 & 0 & -0.1758 & 0.0008 & 0 & -0.6506 & 0 & 0 & 0.0106 & \text{HT-unstable}\\
5 &0 & -0.1776 & 0 & -0.0113 & -0.6473 & 0 & -0.0010 & 0 & \text{HT-unstable}\\
6 & 0 & 0 & -0.1608 & 0.0259 & -0.7994 & 0.1297 & 0 &	0 & \text{HT-unstable}
\end{bmatrix}
$$
where $x_1, x_2, x_3, x_4$ are four coordinates of each \textit{HT-stationary} point and where $g_1, g_2, g_3, g_4$ are the four gradient entries corresponding to each \textit{HT-stationary} point. Since \textit{HT-stationary} points are vectors in $\mathbb{R}^4$, there is no way to show all of them on one 2-d plane. Thus, we use six 2-d plains where each plane shows only two coordinates of \textit{HT-stationary} points. On each 2-d plain we have 6 different points, each one associated with one of the \textit{HT-stationary} points shown in a particular 2-d plain with specified coordinates. In Figure \ref{fig:HT-stationary} the points with red stars are \textit{HT-unstable} ones, and the blue ones are the \textit{HT-stable} ones.
For example, the first 2-d plain (first row-first column) including coordinates $x_1-x_2$ shows the $x_1,x_2$ coordinates of all of the six \textit{HT-stationary} points. On the first row-first column 2-d plain, the first \textit{HT-stationary} point is more distinct because it is the only one that has two nonzeros elements associated with $x_1-x_2$ coordinates. We also can see three points with $x_2=0$, two of which are \textit{HT-unstable} points and one is \textit{HT-unstable} one. This is more clear, if one looks at the column $x_2$ in \textit{HT-stationary} points matrix above. Also, it is clear that we have three \textit{HT-unstable} points with $x_1=0$ on the first 2-d plane.

%
We perturb nonzero coordinates of all \textit{HT-unstable} points with a normal random noise with mean zero and standard deviation of $\sigma=0.5$ to create 4,000 different initialization points. These points create four clouds around \textit{HT-unstable} points which are shown in Figure \ref{fig:unstable_init}.
\begin{figure}[h]
\centering
\includegraphics[scale=0.2]{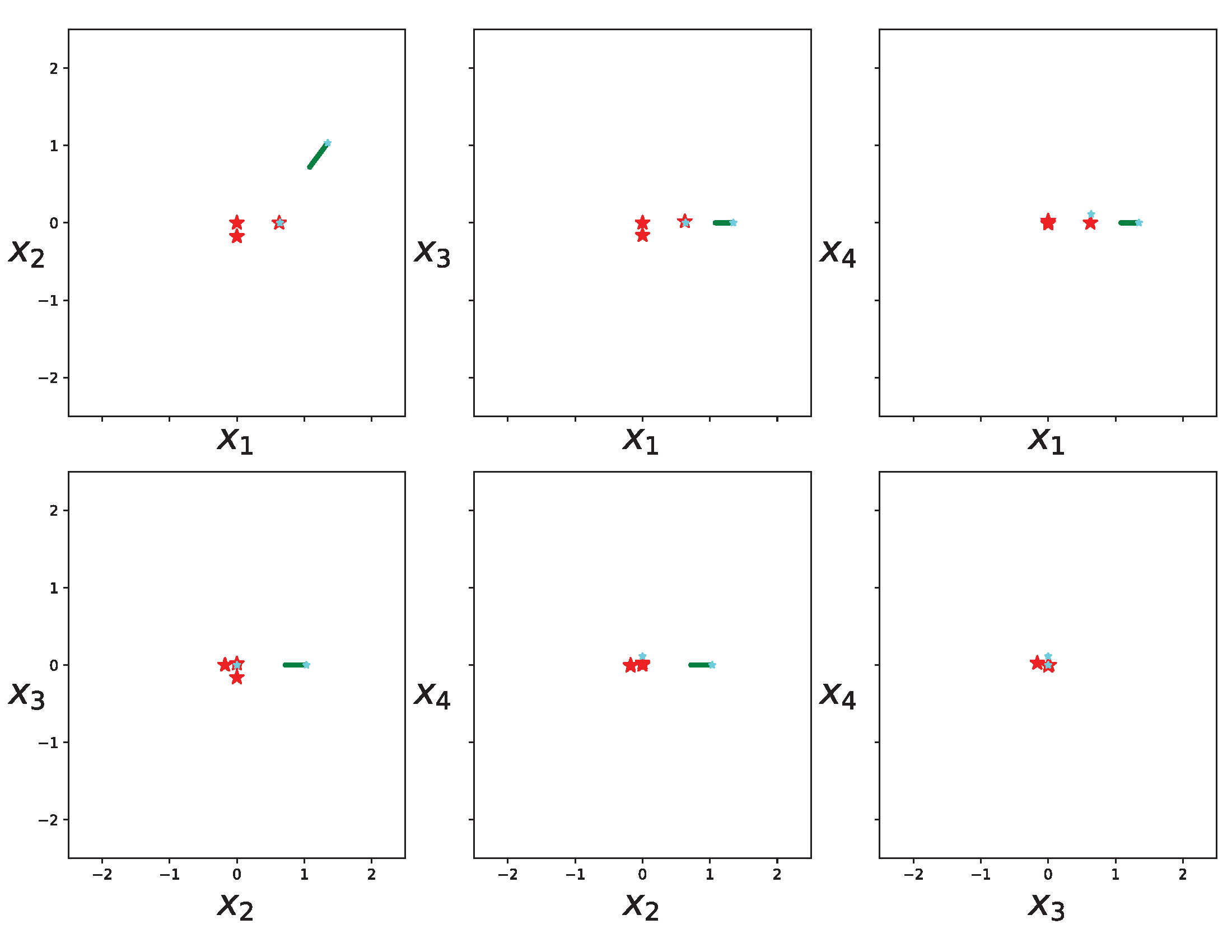}
\caption{Illustration of 4,000 IHT sequences initialized close to four \textit{HT-unstable stationary} points after 400 steps (please refer to the video of all 400 steps). 
}
\label{fig:300thstep}
\end{figure}
Then we run the IHT algorithm for 400 steps.
After 300 steps, all of these initializations escape from those \textit{HT-unstable} points and converge to either of the \textit{HT-stable} stationary points on $x_1-x_2$ or $x_1-x_4$ 2-d planes. In fact, these two \textit{HT-stable} stationary points are sparse local minimizers. Figure \ref{fig:300thstep} shows the 300-th step of IHT algorithm. There is a video in the supplementary materials that shows 400 steps of applying IHT algorithm for 4,000 different runs.
These numerical results corroborate our theoretical results as expected. By looking at the video, one can easily see \textit{escapability} property of \textit{HT-unstable stationary} points, and \textit{Reachability} to \textit{HT-stable stationary} points.

\section{Conclusion}
This paper provide theoretical results that help to understand the IHT algorithm. These theoretical results include a critical gradient descent property of the hard thresholding (HT) operator which is used to show the sequence of the IHT algorithm is decreasing and by doing it over and over we get smaller objective value. This property also allows one to study the IHT algorithm when the stepsize is less than or equal to $1/L_s$, where $L_s>0$ is the Lipschitz constant of the gradient of an objective function. We introduced different stationary points including HT-stable and HT-unstable stationary points and show no matter how close an initialization is to a HT-unstable stationary point, the IHT sequence leaves it. 
We provided a video of 4000 independent runs where the IHT algorithm is initialized very close to a HT-unstable stationary point and showed the sequences escape them.
This property is used to prove that the IHT sequence converges to a HT-stable stationary point. Also, we established a condition for a HT-stable stationary that is a global minimizer with respect to $\gamma=1/L_s$. Finally, we showed the IHT sequence always converges if the function values of HT-stable stationary points are distinct, this is a new assumption that has not been found in the literature.


\newpage
\bibliography{iclr2023_conference}
\bibliographystyle{iclr2023_conference}

\end{document}

%% file: math_commands.tex
\usepackage{amsmath,amsfonts,bm}

\def\eqref#1{equation~\ref{#1}}

\def\1{\bm{1}}

\DeclareMathAlphabet{\mathsfit}{\encodingdefault}{\sfdefault}{m}{sl}
\SetMathAlphabet{\mathsfit}{bold}{\encodingdefault}{\sfdefault}{bx}{n}



%% file: iht.tex
\begin{algorithm}[t]
\caption{The iterative hard thresholding (IHT)}
\label{alg:IHT}
\begin{algorithmic}[1]
\REQUIRE $\mathbf{x}^0\in \mathbb{R}^n$ such that $\|\mathbf{x}^0\|_0\leq s$ and stepsize $\gamma>0$ .
\STATE 
$\mathbf{x}^{k+1} \in H_s(\mathbf{x}^k-\gamma \nabla f(\mathbf{x}^k))$ for $k=0,1, \dots$
\end{algorithmic}
\end{algorithm}

%% file: comparison.tex
\begin{table*}[t]
\centering
\caption{Comparison of results for the deterministic IHT-type algorithms.}
\label{tab:comparison}
\centering
\scalebox{0.6}{
\begin{tabular}{
p{2.4cm}
p{2.8cm}
p{3cm}
p{2.2cm}
p{4.2cm}
p{1cm}
p{7cm}
}
\hline
\toprule
Paper
& Objective function
& Constraints 
& Stepsize 
& Optimality conditions  
& Method
& Convergence
\\
\midrule
\midrule
\cite{beck2013sparsity} 
&  $f(\mathbf{x})$: $L$-LG 
& $\mathbf{x} \in C_s$ 
& $0<\gamma<1/L$ 
& $\mathbf{x}^* \in H_s(\mathbf{x}^* -\frac{1}{L}\nabla f(\mathbf{x}^*))$
& 
IHT
&
\begin{tabular}{@{}c@{}}
Any accumulation point of the IHT sequence
\\
satisfies
$\mathbf{x}^* \in H_s(\mathbf{x}^* -\frac{1}{L}\nabla f(\mathbf{x}^*))$(Theorem 3.1).
\end{tabular}  \\
\midrule
\cite{lu2014iterative}&
\begin{tabular}{@{}c@{}}
$f(\mathbf{x})+\lambda \|\mathbf{x}\|_0$
\\ 
$f(\mathbf{x})$ is convex
\end{tabular}
& 
\begin{tabular}{@{}c@{}}
$\mathbf{x} \in C_s$
\\ 
$l\leq \mathbf{x} \leq u$
\end{tabular}
& $0<\gamma < 1/L$ 
&  Local minimizer &
&
\begin{tabular}{@{}c@{}}
The IHT sequence converges to a local minimizer.
\\ 
(Theorem 3.3)
\end{tabular}
 \\
\midrule
\cite{beck2016minimization}
& $f(\mathbf{x})$: $L$-LG
& \begin{tabular}{@{}c@{}}
$\mathbf{x} \in C_s \cap \mathcal{X}$, $\mathcal{X}$ is 
\\ 
nonnegative symmetric 
\\ or sign free
\end{tabular}
& $0<\gamma < 1/L$ 
&  Basic feasibility
& BFS
& \begin{tabular}{@{}c@{}}
The sequence of BFS converges to a basic 
\\ 
feasible point (Lemma 7.1). 
\end{tabular}

\\
\midrule
\cite{Lu2015OptimizationOS}
& $f(\mathbf{x})$: $L$-LG
& \begin{tabular}{@{}c@{}}
$\mathbf{x} \in C_s \cap \mathcal{X}$, $\mathcal{X}$ is 
\\ 
nonnegative symmetric 
\\ or sign free
\end{tabular}
& $0<\gamma < 1/L$ 
&  $\mathbf{x}^*=P_{C_s \cap \mathcal{X}}(\mathbf{x}^* -\frac{1}{L}\nabla f(\mathbf{x}^*))$ 
& NPG
& \begin{tabular}{@{}c@{}}
An accumulation point of NPG sequence satisfies 
\\
the optimality condition (Theorem 4.3).
\end{tabular} 

\\   
\midrule
\cite{pan2017convergent}
& \begin{tabular}{@{}c@{}}
$f(\mathbf{x})$
\\ 
$L_s$-RSS
\,\&\,
$s$-RC
\end{tabular}
& $\mathbf{x} \in C_s \cap \mathbb{R}^n_{+}$ 
& $0<\gamma < 1/L_s$ 
&  \begin{tabular}{@{}c@{}}
$\mathbf{x}^*\in P_{C_s \cap \mathbb{R}^n_{+}}(\mathbf{x}^* -\frac{1}{L}\nabla f(\mathbf{x}^*))$
\\ 
or $L_s$-stationary
\end{tabular}
&
IIHT
& \begin{tabular}{@{}c@{}}
Any accumulation point of IIHT sequence \\
converges to a $L_s$-stationary point (Theorem 3.1).
\end{tabular} 
\\   
\midrule
\cite{pan2017convergent}
& \begin{tabular}{@{}c@{}}
$f(\mathbf{x})$
\\ 
$L_s$-RSS
\,\&\,
$\beta_s$-RSC
\end{tabular}
& $\mathbf{x} \in C_s \cap \mathbb{R}^n_{+}$ & $0<\gamma < 1/L_s$ 
&  \begin{tabular}{@{}c@{}}
$\mathbf{x}^*\in P_{C_s \cap \mathbb{R}^n_{+}}(\mathbf{x}^* -\frac{1}{L}\nabla f(\mathbf{x}^*))$
\\ 
or $L_s$-stationary
\end{tabular}
&
IIHT
& \begin{tabular}{@{}c@{}}
The IIHT sequence converges to  a local \\
minimizer (Theorem 3.2). If $\|\mathbf{x}^*\|_0<s$, then \\
$\mathbf{x}^*$ is a global minimizer. When $\|\mathbf{x}^*\|_0=s$,
\\
the IIHT sequence converges Q-linearly 
\\
(Theorem 3.4).
\end{tabular} 
\\
\midrule
\cite{zhou2021global}
& \begin{tabular}{@{}c@{}}
$f(\mathbf{x})$
\\ 
$L_s$-RSS
\end{tabular}
& $\mathbf{x} \in C_s$ 
& $0<\gamma < 1/L_s$ 
&  \begin{tabular}{@{}c@{}}
$\mathbf{x}^*\in H_s(\mathbf{x}^* -\frac{1}{L_s}\nabla f(\mathbf{x}^*))$ 
\\ 
or $L_s$-stationary
\end{tabular}
&
NHTP
&  \begin{tabular}{@{}c@{}}
Any accumulation point $\mathbf{x}^*$ of the NHTP  \\
sequence is an $L_s$-stationary point (Theorem 9). 
\\ 
If $\mathbf{x}^*$ is isolated, 
the entire sequence converges.
\end{tabular}
\\
\midrule
\cite{zhou2021global}
& \begin{tabular}{@{}c@{}}
$f(\mathbf{x})$
\\ 
$L_s$-RSS
\,\&\,
$\beta_s$-RSC
\\
\& Restricted Hessian
\\
is Lipschitz
\end{tabular}
& $\mathbf{x} \in C_s$ 
& $0<\gamma < 1/L_s$ 
& \begin{tabular}{@{}c@{}}
$\mathbf{x}^*\in H_s(\mathbf{x}^* -\frac{1}{L_s}\nabla f(\mathbf{x}^*))$ 
\\ 
or $L_s$-stationary
\end{tabular} 
&
NHTP
&  \begin{tabular}{@{}c@{}}
The NHTP sequence converges to a $L_s$-stationary
\\ point (Theorem 10).
\\
Locally, it converges quadratically.
\end{tabular}
\\
   
\midrule
Ours
& \begin{tabular}{@{}c@{}}
$f(\mathbf{x})$
\\ 
$L_s$-RSS
\end{tabular}
& $\mathbf{x} \in C_s$ 
& $0<\gamma \leq 1/L_s$ 
& \begin{tabular}{@{}c@{}}
$\mathbf{x}^*\in H_s(\mathbf{x}^* -\frac{1}{L_s}\nabla f(\mathbf{x}^*))$ 
\\ 
or HT-stable
\end{tabular}
&
IHT
& \begin{tabular}{@{}l@{}}
If $f(\mathbf{x}^*) \neq f(\mathbf{y}^*)$ for all $\mathbf{x}^*,\mathbf{y}^*$ HT-stable points,
\\
by starting from a $\mathbf{x}^0 \in C_s$ the IHT sequence 
\\ 
converges to some HT-stable points. 
\\
Corollary \ref{cor:IHTconvergence}: global convergence
\end{tabular}
\\
\midrule
Ours
& \begin{tabular}{@{}c@{}}
$f(\mathbf{x})$
\\ 
$L_s$-RSS
\end{tabular}
& $\mathbf{x} \in C_s$ 
& $0<\gamma < 1/L_s$ 
& \begin{tabular}{@{}c@{}}
$\mathbf{x}^*\in H_s(\mathbf{x}^* -\frac{1}{L_s}\nabla f(\mathbf{x}^*))$ 
\\ 
or HT-stable
\end{tabular}
&
IHT
& \begin{tabular}{@{}l@{}}
If all HT-stable points are isolated,
\\
by starting from a $\mathbf{x}^0 \in C_s$ the IHT sequence 
\\ 
converges to some HT-stable points. 
\\
Corollary \ref{cor:IHTconvergence}: global convergence
\end{tabular}
\\
\midrule
Ours
& \begin{tabular}{@{}c@{}}
$f(\mathbf{x})$
\\ 
$L_s$-RSS 
\\
Strictly convex on $\mathcal{S_{\mathcal{J}}}$
\\
with $|\mathcal{J}|=s$
\end{tabular}
& $\mathbf{x} \in C_s$ 
& $0<\gamma \leq 1/L_s$ 
& \begin{tabular}{@{}c@{}}
$\mathbf{x}^*= H_s(\mathbf{x}^* -\frac{1}{L_s}\nabla f(\mathbf{x}^*))$ 
\\ 
or HT-strictly stable
\end{tabular}
&
IHT
& \begin{tabular}{@{}l@{}}
By starting from $\mathbf{x}^0 \in \mathcal{B}(\mathbf{x}^*, \delta) \cap C_s$ the IHT 
\\
sequence converges to HT-strictly stable $\mathbf{x}^*$ point
\\ 
Q-linearly (Proposition \ref{prop:reachability}).
\end{tabular}
\\
\bottomrule
\hline
\\
\end{tabular}
}
\end{table*}

%% file: iclr2023_conference.bbl
\begin{thebibliography}{64}
\providecommand{\natexlab}[1]{#1}
\providecommand{\url}[1]{\texttt{#1}}
\expandafter\ifx\csname urlstyle\endcsname\relax
  \providecommand{\doi}[1]{doi: #1}\else
  \providecommand{\doi}{doi: \begingroup \urlstyle{rm}\Url}\fi

\bibitem[Bahmani et~al.(2013)Bahmani, Raj, and Boufounos]{bahmani2013greedy}
Sohail Bahmani, Bhiksha Raj, and Petros~T Boufounos.
\newblock Greedy sparsity-constrained optimization.
\newblock \emph{Journal of Machine Learning Research}, 14\penalty0
  (Mar):\penalty0 807--841, 2013.

\bibitem[Bauschke et~al.(2014)Bauschke, Luke, Phan, and
  Wang]{bauschke2014restricted}
Heinz~H Bauschke, D~Russell Luke, Hung~M Phan, and Xianfu Wang.
\newblock Restricted normal cones and sparsity optimization with affine
  constraints.
\newblock \emph{Foundations of Computational Mathematics}, 14\penalty0
  (1):\penalty0 63--83, 2014.

\bibitem[Beck \& Eldar(2013)Beck and Eldar]{beck2013sparsity}
Amir Beck and Yonina~C Eldar.
\newblock Sparsity constrained nonlinear optimization: Optimality conditions
  and algorithms.
\newblock \emph{SIAM Journal on Optimization}, 23\penalty0 (3):\penalty0
  1480--1509, 2013.

\bibitem[Beck \& Hallak(2016)Beck and Hallak]{beck2016minimization}
Amir Beck and Nadav Hallak.
\newblock On the minimization over sparse symmetric sets: projections,
  optimality conditions, and algorithms.
\newblock \emph{Mathematics of Operations Research}, 41\penalty0 (1):\penalty0
  196--223, 2016.

\bibitem[Blumensath(2012)]{blumensath2012accelerated}
Thomas Blumensath.
\newblock Accelerated iterative hard thresholding.
\newblock \emph{Signal Processing}, 92\penalty0 (3):\penalty0 752--756, 2012.

\bibitem[Blumensath \& Davies(2008)Blumensath and
  Davies]{blumensath2008iterative}
Thomas Blumensath and Mike~E Davies.
\newblock Iterative thresholding for sparse approximations.
\newblock \emph{Journal of Fourier analysis and Applications}, 14\penalty0
  (5):\penalty0 629--654, 2008.

\bibitem[Blumensath \& Davies(2009)Blumensath and
  Davies]{blumensath2009iterative}
Thomas Blumensath and Mike~E Davies.
\newblock Iterative hard thresholding for compressed sensing.
\newblock \emph{Applied and computational harmonic analysis}, 27\penalty0
  (3):\penalty0 265--274, 2009.

\bibitem[Bourguignon et~al.(2015)Bourguignon, Ninin, Carfantan, and
  Mongeau]{bourguignon2015exact}
S{\'e}bastien Bourguignon, Jordan Ninin, Herv{\'e} Carfantan, and Marcel
  Mongeau.
\newblock Exact sparse approximation problems via mixed-integer programming:
  Formulations and computational performance.
\newblock \emph{IEEE Transactions on Signal Processing}, 64\penalty0
  (6):\penalty0 1405--1419, 2015.

\bibitem[Brodie et~al.(2009)Brodie, Daubechies, De~Mol, Giannone, and
  Loris]{brodie2009sparse}
Joshua Brodie, Ingrid Daubechies, Christine De~Mol, Domenico Giannone, and
  Ignace Loris.
\newblock Sparse and stable markowitz portfolios.
\newblock \emph{Proceedings of the National Academy of Sciences}, 106\penalty0
  (30):\penalty0 12267--12272, 2009.

\bibitem[Bucher \& Schwartz(2018)Bucher and Schwartz]{bucher2018second}
Max Bucher and Alexandra Schwartz.
\newblock Second-order optimality conditions and improved convergence results
  for regularization methods for cardinality-constrained optimization problems.
\newblock \emph{Journal of Optimization Theory and Applications}, 178\penalty0
  (2):\penalty0 383--410, 2018.

\bibitem[Burdakov et~al.(2016)Burdakov, Kanzow, and
  Schwartz]{burdakov2016mathematical}
Oleg~P Burdakov, Christian Kanzow, and Alexandra Schwartz.
\newblock Mathematical programs with cardinality constraints: reformulation by
  complementarity-type conditions and a regularization method.
\newblock \emph{SIAM Journal on Optimization}, 26\penalty0 (1):\penalty0
  397--425, 2016.

\bibitem[Chang et~al.(2000)Chang, Meade, Beasley, and
  Sharaiha]{chang2000heuristics}
T-J Chang, Nigel Meade, John~E Beasley, and Yazid~M Sharaiha.
\newblock Heuristics for cardinality constrained portfolio optimisation.
\newblock \emph{Computers \& Operations Research}, 27\penalty0 (13):\penalty0
  1271--1302, 2000.

\bibitem[Chartrand(2007)]{chartrand2007exact}
Rick Chartrand.
\newblock Exact reconstruction of sparse signals via nonconvex minimization.
\newblock \emph{IEEE Signal Processing Letters}, 14\penalty0 (10):\penalty0
  707--710, 2007.

\bibitem[Chen \& Gu(2016)Chen and Gu]{chen2016accelerated}
Jinghui Chen and Quanquan Gu.
\newblock Accelerated stochastic block coordinate gradient descent for sparsity
  constrained nonconvex optimization.
\newblock In \emph{UAI}, 2016.

\bibitem[Chun \& Kele{\c{s}}(2010)Chun and Kele{\c{s}}]{chun2010sparse}
Hyonho Chun and S{\"u}nd{\"u}z Kele{\c{s}}.
\newblock Sparse partial least squares regression for simultaneous dimension
  reduction and variable selection.
\newblock \emph{Journal of the Royal Statistical Society: Series B (Statistical
  Methodology)}, 72\penalty0 (1):\penalty0 3--25, 2010.

\bibitem[Damadi et~al.(2022)Damadi, Nouri, and Pirsiavash]{damadi2022amenable}
Saeed Damadi, Erfan Nouri, and Hamed Pirsiavash.
\newblock Amenable sparse network investigator.
\newblock \emph{arXiv preprint arXiv:2202.09284}, 2022.

\bibitem[Davenport \& Romberg(2016)Davenport and
  Romberg]{davenport2016overview}
Mark~A Davenport and Justin Romberg.
\newblock An overview of low-rank matrix recovery from incomplete observations.
\newblock \emph{IEEE Journal of Selected Topics in Signal Processing},
  10\penalty0 (4):\penalty0 608--622, 2016.

\bibitem[Davis(1994)]{davis1994adaptive}
Geoffrey~Mark Davis.
\newblock \emph{Adaptive nonlinear approximations}.
\newblock PhD thesis, New York University, 1994.

\bibitem[Dedieu et~al.(2021)Dedieu, Hazimeh, and Mazumder]{dedieu2021learning}
Antoine Dedieu, Hussein Hazimeh, and Rahul Mazumder.
\newblock Learning sparse classifiers: Continuous and mixed integer
  optimization perspectives.
\newblock \emph{Journal of Machine Learning Research}, 22\penalty0
  (135):\penalty0 1--47, 2021.

\bibitem[Desboulets(2018)]{desboulets2018review}
Loann David~Denis Desboulets.
\newblock A review on variable selection in regression analysis.
\newblock \emph{Econometrics}, 6\penalty0 (4):\penalty0 45, 2018.

\bibitem[Donoho et~al.(2012)Donoho, Tsaig, Drori, and Starck]{donoho2012sparse}
David~L Donoho, Yaakov Tsaig, Iddo Drori, and Jean-Luc Starck.
\newblock Sparse solution of underdetermined systems of linear equations by
  stagewise orthogonal matching pursuit.
\newblock \emph{IEEE transactions on Information Theory}, 58\penalty0
  (2):\penalty0 1094--1121, 2012.

\bibitem[Eldar \& Kutyniok(2012)Eldar and Kutyniok]{eldar2012compressed}
Yonina~C Eldar and Gitta Kutyniok.
\newblock \emph{Compressed sensing: theory and applications}.
\newblock Cambridge university press, 2012.

\bibitem[Fan \& Li(2001)Fan and Li]{fan2001variable}
Jianqing Fan and Runze Li.
\newblock Variable selection via nonconcave penalized likelihood and its oracle
  properties.
\newblock \emph{Journal of the American statistical Association}, 96\penalty0
  (456):\penalty0 1348--1360, 2001.

\bibitem[Fornasier et~al.(2011)Fornasier, Rauhut, and Ward]{fornasier2011low}
Massimo Fornasier, Holger Rauhut, and Rachel Ward.
\newblock Low-rank matrix recovery via iteratively reweighted least squares
  minimization.
\newblock \emph{SIAM Journal on Optimization}, 21\penalty0 (4):\penalty0
  1614--1640, 2011.

\bibitem[Foucart \& Lai(2009)Foucart and Lai]{foucart2009sparsest}
Simon Foucart and Ming-Jun Lai.
\newblock Sparsest solutions of underdetermined linear systems via
  $\ell_q$-minimization for $0< q\leq 1$.
\newblock \emph{Applied and Computational Harmonic Analysis}, 26\penalty0
  (3):\penalty0 395--407, 2009.

\bibitem[Foucart \& Rauhut(2013)Foucart and Rauhut]{foucart2mathematical}
Simon Foucart and Holger Rauhut.
\newblock A mathematical introduction to compressive sensing.
\newblock pp.\  1--39. Springer, 2013.

\bibitem[Gale et~al.(2019)Gale, Elsen, and Hooker]{gale2019state}
Trevor Gale, Erich Elsen, and Sara Hooker.
\newblock The state of sparsity in deep neural networks.
\newblock \emph{arXiv preprint arXiv:1902.09574}, 2019.

\bibitem[Haeffele et~al.(2014)Haeffele, Young, and
  Vidal]{haeffele2014structured}
Benjamin Haeffele, Eric Young, and Rene Vidal.
\newblock Structured low-rank matrix factorization: Optimality, algorithm, and
  applications to image processing.
\newblock In \emph{International conference on machine learning}, pp.\
  2007--2015. PMLR, 2014.

\bibitem[Jain et~al.(2014)Jain, Tewari, and Kar]{jain2014iterative}
Prateek Jain, Ambuj Tewari, and Purushottam Kar.
\newblock On iterative hard thresholding methods for high-dimensional
  m-estimation.
\newblock \emph{Advances in neural information processing systems}, 27, 2014.

\bibitem[Khanna \& Kyrillidis(2018)Khanna and Kyrillidis]{khanna2018iht}
Rajiv Khanna and Anastasios Kyrillidis.
\newblock Iht dies hard: Provable accelerated iterative hard thresholding.
\newblock In \emph{International Conference on Artificial Intelligence and
  Statistics}, pp.\  188--198. PMLR, 2018.

\bibitem[Lai \& Wang(2011)Lai and Wang]{lai2011unconstrained}
Ming-Jun Lai and Jingyue Wang.
\newblock An unconstrained $\ell_q$ minimization with $0<q \leq< 1$ for sparse
  solution of underdetermined linear systems.
\newblock \emph{SIAM Journal on Optimization}, 21\penalty0 (1):\penalty0
  82--101, 2011.

\bibitem[Li et~al.(2016)Li, Arora, Liu, Haupt, and Zhao]{li2016nonconvex}
Xingguo Li, Raman Arora, Han Liu, Jarvis Haupt, and Tuo Zhao.
\newblock Nonconvex sparse learning via stochastic optimization with
  progressive variance reduction.
\newblock \emph{arXiv preprint arXiv:1605.02711}, 2016.

\bibitem[Li \& Song(2015)Li and Song]{li2015first}
Xue Li and Wen Song.
\newblock The first-order necessary conditions for sparsity constrained
  optimization.
\newblock \emph{Journal of the Operations Research Society of China},
  3\penalty0 (4):\penalty0 521--535, 2015.

\bibitem[Liang et~al.(2020)Liang, Tong, Zhu, and Bi]{liang2020effective}
Guannan Liang, Qianqian Tong, Chunjiang Zhu, and Jinbo Bi.
\newblock An effective hard thresholding method based on stochastic variance
  reduction for nonconvex sparse learning.
\newblock In \emph{Proceedings of the AAAI Conference on Artificial
  Intelligence}, volume~34, pp.\  1585--1592, 2020.

\bibitem[Liu et~al.(2017)Liu, Yuan, Wang, Liu, and Metaxas]{liu2017dual}
Bo~Liu, Xiao-Tong Yuan, Lezi Wang, Qingshan Liu, and Dimitris~N Metaxas.
\newblock Dual iterative hard thresholding: From non-convex sparse minimization
  to non-smooth concave maximization.
\newblock In \emph{International Conference on Machine Learning}, pp.\
  2179--2187. PMLR, 2017.

\bibitem[Liu \& Foygel~Barber(2020)Liu and Foygel~Barber]{liu2020between}
Haoyang Liu and Rina Foygel~Barber.
\newblock Between hard and soft thresholding: optimal iterative thresholding
  algorithms.
\newblock \emph{Information and Inference: A Journal of the IMA}, 9\penalty0
  (4):\penalty0 899--933, 2020.

\bibitem[Lu(2014)]{lu2014iterative}
Zhaosong Lu.
\newblock Iterative hard thresholding methods for $l_0$ regularized convex cone
  programming.
\newblock \emph{Mathematical Programming}, 147\penalty0 (1):\penalty0 125--154,
  2014.

\bibitem[Lu(2015)]{Lu2015OptimizationOS}
Zhaosong Lu.
\newblock Optimization over sparse symmetric sets via a nonmonotone projected
  gradient method.
\newblock \emph{arXiv: Optimization and Control}, 2015.

\bibitem[Mallat \& Zhang(1993)Mallat and Zhang]{mallat1993matching}
St{\'e}phane~G Mallat and Zhifeng Zhang.
\newblock Matching pursuits with time-frequency dictionaries.
\newblock \emph{IEEE Transactions on signal processing}, 41\penalty0
  (12):\penalty0 3397--3415, 1993.

\bibitem[Molchanov et~al.(2017)Molchanov, Ashukha, and
  Vetrov]{molchanov2017variational}
Dmitry Molchanov, Arsenii Ashukha, and Dmitry Vetrov.
\newblock Variational dropout sparsifies deep neural networks.
\newblock In \emph{International Conference on Machine Learning}, pp.\
  2498--2507. PMLR, 2017.

\bibitem[Mousavi \& Shen(2019)Mousavi and Shen]{mousavi2019solution}
Seyedahmad Mousavi and Jinglai Shen.
\newblock Solution uniqueness of convex piecewise affine functions based
  optimization with applications to constrained $\ell_1$ minimization.
\newblock \emph{ESAIM: Control, Optimisation and Calculus of Variations},
  25:\penalty0 56, 2019.

\bibitem[Natarajan(1995)]{natarajan1995sparse}
Balas~Kausik Natarajan.
\newblock Sparse approximate solutions to linear systems.
\newblock \emph{SIAM journal on computing}, 24\penalty0 (2):\penalty0 227--234,
  1995.

\bibitem[Needell \& Tropp(2009)Needell and Tropp]{needell2009cosamp}
Deanna Needell and Joel~A Tropp.
\newblock Cosamp: Iterative signal recovery from incomplete and inaccurate
  samples.
\newblock \emph{Applied and computational harmonic analysis}, 26\penalty0
  (3):\penalty0 301--321, 2009.

\bibitem[Needell \& Vershynin(2009)Needell and Vershynin]{needell2009uniform}
Deanna Needell and Roman Vershynin.
\newblock Uniform uncertainty principle and signal recovery via regularized
  orthogonal matching pursuit.
\newblock \emph{Foundations of computational mathematics}, 9\penalty0
  (3):\penalty0 317--334, 2009.

\bibitem[Needell \& Vershynin(2010)Needell and Vershynin]{needell2010signal}
Deanna Needell and Roman Vershynin.
\newblock Signal recovery from incomplete and inaccurate measurements via
  regularized orthogonal matching pursuit.
\newblock \emph{IEEE Journal of selected topics in signal processing},
  4\penalty0 (2):\penalty0 310--316, 2010.

\bibitem[Negahban et~al.(2012)Negahban, Ravikumar, Wainwright, and
  Yu]{negahban2012unified}
Sahand~N Negahban, Pradeep Ravikumar, Martin~J Wainwright, and Bin Yu.
\newblock A unified framework for high-dimensional analysis of $ m $-estimators
  with decomposable regularizers.
\newblock \emph{Statistical science}, 27\penalty0 (4):\penalty0 538--557, 2012.

\bibitem[Pan et~al.(2015)Pan, Xiu, and Zhou]{pan2015solutions}
Li-Li Pan, Nai-Hua Xiu, and Sheng-Long Zhou.
\newblock On solutions of sparsity constrained optimization.
\newblock \emph{Journal of the Operations Research Society of China},
  3\penalty0 (4):\penalty0 421--439, 2015.

\bibitem[Pan et~al.(2017{\natexlab{a}})Pan, Xiu, and Fan]{pan2017optimality}
LiLi Pan, NaiHua Xiu, and Jun Fan.
\newblock Optimality conditions for sparse nonlinear programming.
\newblock \emph{Science China Mathematics}, 60\penalty0 (5):\penalty0 759--776,
  2017{\natexlab{a}}.

\bibitem[Pan et~al.(2017{\natexlab{b}})Pan, Zhou, Xiu, and
  Qi]{pan2017convergent}
Lili Pan, Shenglong Zhou, Naihua Xiu, and Houduo Qi.
\newblock A convergent iterative hard thresholding for sparsity and
  nonnegativity constrained optimization.
\newblock \emph{Pacific Journal of Optimization}, 13\penalty0 (2):\penalty0
  325--353, 2017{\natexlab{b}}.

\bibitem[Paszke et~al.(2019)Paszke, Gross, Massa, Lerer, Bradbury, Chanan,
  Killeen, Lin, Gimelshein, Antiga, et~al.]{paszke2019pytorch}
Adam Paszke, Sam Gross, Francisco Massa, Adam Lerer, James Bradbury, Gregory
  Chanan, Trevor Killeen, Zeming Lin, Natalia Gimelshein, Luca Antiga, et~al.
\newblock Pytorch: An imperative style, high-performance deep learning library.
\newblock \emph{arXiv preprint arXiv:1912.01703}, 2019.

\bibitem[Pati et~al.(1993)Pati, Rezaiifar, and
  Krishnaprasad]{pati1993orthogonal}
Yagyensh~Chandra Pati, Ramin Rezaiifar, and Perinkulam~Sambamurthy
  Krishnaprasad.
\newblock Orthogonal matching pursuit: Recursive function approximation with
  applications to wavelet decomposition.
\newblock In \emph{Proceedings of 27th Asilomar conference on signals, systems
  and computers}, pp.\  40--44. IEEE, 1993.

\bibitem[Tibshirani(1996)]{tibshirani1996regression}
Robert Tibshirani.
\newblock Regression shrinkage and selection via the lasso.
\newblock \emph{Journal of the Royal Statistical Society: Series B
  (Methodological)}, 58\penalty0 (1):\penalty0 267--288, 1996.

\bibitem[Vu \& Raich(2019)Vu and Raich]{vu2019accelerating}
Trung Vu and Raviv Raich.
\newblock Accelerating iterative hard thresholding for low-rank matrix
  completion via adaptive restart.
\newblock In \emph{ICASSP 2019-2019 IEEE International Conference on Acoustics,
  Speech and Signal Processing (ICASSP)}, pp.\  2917--2921. IEEE, 2019.

\bibitem[Wang et~al.(2011)Wang, Xu, and Tang]{wang2011performance}
Meng Wang, Weiyu Xu, and Ao~Tang.
\newblock On the performance of sparse recovery via $\ell_p$-minimization $0<p
  \leq 1$.
\newblock \emph{IEEE Transactions on Information Theory}, 57\penalty0
  (11):\penalty0 7255--7278, 2011.

\bibitem[Wang et~al.(2019)Wang, Xiu, and Zhou]{wang2019fast}
Rui Wang, Naihua Xiu, and Shenglong Zhou.
\newblock Fast newton method for sparse logistic regression.
\newblock \emph{arXiv}, \penalty0 (1901.02768), 2019.

\bibitem[Won et~al.(2022)Won, Lange, and Xu]{won2022unified}
Joong-Ho Won, Kenneth Lange, and Jason Xu.
\newblock A unified analysis of convex and non-convex lp-ball projection
  problems.
\newblock \emph{arXiv preprint arXiv:2203.00564}, 2022.

\bibitem[Wu \& Bian(2020)Wu and Bian]{wu2020accelerated}
Fan Wu and Wei Bian.
\newblock Accelerated iterative hard thresholding algorithm for $l_0$
  regularized regression problem.
\newblock \emph{Journal of Global Optimization}, 76\penalty0 (4):\penalty0
  819--840, 2020.

\bibitem[Zhao et~al.(2021)Zhao, Xiu, Qi, and Luo]{zhao2021lagrange}
Chen Zhao, Naihua Xiu, Houduo Qi, and Ziyan Luo.
\newblock A lagrange--newton algorithm for sparse nonlinear programming.
\newblock \emph{Mathematical Programming}, pp.\  1--26, 2021.

\bibitem[Zheng et~al.(2017)Zheng, Maleki, Weng, Wang, and Long]{zheng2017does}
Le~Zheng, Arian Maleki, Haolei Weng, Xiaodong Wang, and Teng Long.
\newblock Does $\ell_p$-minimization outperform $\ell_1$-minimization?
\newblock \emph{IEEE Transactions on Information Theory}, 63\penalty0
  (11):\penalty0 6896--6935, 2017.

\bibitem[Zhou et~al.(2018)Zhou, Yuan, and Feng]{zhou2018efficient}
Pan Zhou, Xiaotong Yuan, and Jiashi Feng.
\newblock Efficient stochastic gradient hard thresholding.
\newblock \emph{Advances in Neural Information Processing Systems}, 31, 2018.

\bibitem[Zhou et~al.(2020)Zhou, Pan, and Xiu]{zhou2020subspace}
Shenglong Zhou, Lili Pan, and Naihua Xiu.
\newblock Subspace newton method for the $l_0$-regularized optimization.
\newblock \emph{arXiv}, \penalty0 (2004.05132), 2020.

\bibitem[Zhou et~al.(2021)Zhou, Xiu, and Qi]{zhou2021global}
Shenglong Zhou, Naihua Xiu, and Hou-Duo Qi.
\newblock Global and quadratic convergence of newton hard-thresholding pursuit.
\newblock \emph{J. Mach. Learn. Res.}, 22\penalty0 (12):\penalty0 1--45, 2021.

\bibitem[Zhu et~al.(2018)Zhu, Dong, Yu, and Chen]{zhu2018lagrange}
Wenxing Zhu, Zhengshan Dong, Yuanlong Yu, and Jianli Chen.
\newblock Lagrange dual method for sparsity constrained optimization.
\newblock \emph{IEEE Access}, 6:\penalty0 28404--28416, 2018.

\bibitem[Zou \& Hastie(2005)Zou and Hastie]{zou2005regularization}
Hui Zou and Trevor Hastie.
\newblock Regularization and variable selection via the elastic net.
\newblock \emph{Journal of the royal statistical society: series B (statistical
  methodology)}, 67\penalty0 (2):\penalty0 301--320, 2005.

\end{thebibliography}
